\documentclass[10pt,twoside]{article}
\usepackage{amssymb,cite}
\usepackage{graphicx}
\usepackage{amsmath}
\usepackage{Latex-document}

\newcommand{\Rot}{\mbox{SO}(3)}
\newcommand{\Mat}{{\mathbb R}^{3 \times 3}}
\newcommand{\Vector}{{\mathbb R}^3}
\newcommand{\A}{{\mathcal A}}
\newcommand{\U}{{\mathcal U}}
\newcommand{\V}{{\mathcal V}}
\newcommand{\gy}{\nabla y}
\newcommand{\meas}{\operatorname{meas}}
\newcommand{\esssup}{\operatornamewithlimits{ess\,sup}}
\newcommand{\E}{{\mathcal E}}
\newcommand{\Real}{{\mathbb R}}
\def\jump#1{\lbrack\!\lbrack\,#1\,\rbrack\!\rbrack}

\title{\bf Computational Modeling of\vskip -2mm Microstructure\thanks{This
work was supported in part by AFOSR F49620-98-1-0433, by NSF
DMS-0074043, by ARO DAAG55-98-1-0335, by the Caltech CIMMS, and by
the Minnesota Supercomputer Institute.}\vskip 6mm}

\author{Mitchell Luskin\vspace*{-0.5cm}\thanks{School of Mathematics,
University of Minnesota, Minneapolis, Minnesota 55455, USA.
E-mail: luskin@math.umn.edu}}
\date{\vspace{-8mm}}

\begin{document}

\maketitle

\thispagestyle{first} \setcounter{page}{707}

\begin{abstract}\vskip 3mm
Many materials such as martensitic
or ferromagnetic crystals
are observed to be in metastable states exhibiting
a fine-scale, structured spatial oscillation called microstructure;
and hysteresis is observed as the temperature,
boundary forces, or external magnetic field changes.
We have developed a numerical analysis of microstructure and
used this theory to construct numerical methods
that have been used to compute approximations to
the deformation of crystals with microstructure.
\vskip 4.5mm


\end{abstract}

\vskip 12mm

\section{Introduction} \label{section 1}\setzero
\vskip-5mm \hspace{5mm }

Martensitic crystals are observed to be in metastable
states that can be modeled
by local minima of the energy
[1,\,2,\,11,\,17,\,19,\,25,\,33,\,36]
\begin{equation}
  \mathcal E(y) = \int_{\Omega} \phi(\nabla y(x),\theta(x)) \, dx
+\mbox{interfacial energy}+\mbox{loading energy}, \label{energy}
\end{equation}
where $\Omega\subset{\mathbb R}^3$
 is the reference configuration of the crystal,
$y(x):\Omega\to\mathbb R^3$ is the deformation that may be constrained on
the boundary $\partial\Omega,$ and
$\theta(x):\Omega\to\mathbb R$ is the temperature.
The frame-indifferent elastic energy density
$\phi(F,\theta):\mathbb R^{3 \times 3}\times \mathbb R\to
\mathbb R$ is minimized at high temperature
$\theta\ge\theta_T$ on $\Rot$
and at low temperature $\theta\le\theta_T$ on the martensitic variants
$\mathcal U=\Rot U_1\cup\dots\cup\Rot U_N$ where
the $U_i\in\mathbb R^{3 \times 3}$ are symmetry-related
transformation strains satisfying
\[
  \{ \,R_i^T U_1 R_i: R_i \in \mathcal G \,\} = \{\, U_{1}, \, \ldots, U_{N}\,\}
\]
for the symmetry group $\mathcal G$ of the high temperature (austenitic) phase.
The loading energy above results from applied boundary
forces.

Microstructure occurs when the deformation
gradient oscillates in space among the
$\Rot U_i$ to enable the deformation
to attain a lower energy
than could be attained by a more homogeneous state [1,\,25].
The simplest microstructure is a laminate in which
the deformation gradient oscillates between
$R_iU_i\in\Rot U_i$ and $R_jU_j\in\Rot U_j$ for $i\ne j$ in parallel layers
of fine scale, but more complex microstructure
is observed in nature and is predicted by the
theory [2,\,25].

We have developed numerical methods for the
computation of microstructure in martensitic
and ferromagnetic crystals and validated these
methods by the development of a numerical analysis of
microstructure
[4,\,6,\,12,\,14,\,16,\,22,\,25--28].
Related results are given in
[9,\,10,\,15,\,21,\,24,\,31,\,32,\,34].
For
martensitic crystals, we have given error estimates for stable quantities such as
nonlinear integrals
$\int_\Omega f(x,\nabla y(x))\,dx$ for smooth
functions $f(x,F):\Omega\times \mathbb R^{3 \times 3}
\to\mathbb R$ and for the local volume fractions
(Young measure) of the variants $\Rot U_i$
even though pointwise values of the deformation
gradient are not stable under mesh refinement.

To model the evolution of metastable states,
we have developed a computational model that nucleates the first order phase change
since otherwise the
crystal would remain stuck in local minima
of the energy as the temperature
or boundary forces are varied [8].  Our finite element
model for the quasi-static evolution of the
martensitic phase transformation in a thin film
nucleates regions of the high temperature
phase during heating and regions of the low temperature
phase during cooling.

Graphical images for the computations of microstructure and phase transformation
described in this paper can be found at http://www.math.umn.edu/\verb+~+luskin
and in the cited references.
A more extensive description of microstructure and
its computation can be found at the above website
as well as in the selected references at the end of this paper.

\section{Numerical analysis of microstructure} \label{section 2}\setzero
\vskip-5mm \hspace{5mm }

Martensitic crystals typically exist in metastable states for
time-scales of technological interest.  Many important analytic
results have been obtained for mathematical models for martensitic crystals,
especially for energy-minimizing deformations with microstructure
[1,\,2,\,11,\,20,\,30,\,35].
These results and concepts for energy-minimizing deformations
should also have a role in the analysis of metastability
[18,\,29].  Similarly,
we have developed a numerical analysis of
microstructure
[4,\,6,\,12,\,14,\,16,\,22,\,25--28]
for which results have been obtained primarily for the approximation of
energy-minimizing deformations that we think also give insight and some
validation for the investigation of metastability by computational methods.

We give here a summary of the numerical analysis of martensitic microstructure
that we have developed for temperatures $\theta<\theta_T$ for which the energy density
$\phi(F,\theta)$ is minimized on the martensitic
variants $\mathcal U=\Rot U_1\cup\dots\cup\Rot U_N.$

We  assume that the energy density $\phi(F,\theta)$ is continuous and satisfies
near the minimizing deformation gradients
$\U$ the
quadratic growth condition  given by
\begin{equation}
\label{eq:growth_condition}
  \phi(F,\theta) \ge \mu \, \| F - \pi(F) \|^{2} \qquad \text{for all } F \in
  \Mat,
\end{equation}
where $\mu > 0$ is a constant and $\pi: \Mat \rightarrow \U$ is a
projection satisfying
\begin{equation*}
\label{eq:definition_pi}
  \| F - \pi(F) \| = \min_{G \in \U} \| F - G \| \qquad \text{for all }
  F \in \Mat.
\end{equation*}
We also assume that the energy density
$\phi(F,\theta)$ satisfies the growth condition for large $F$ given by
\[
\phi(F,\theta)\ge C_1\|F\|^p-C_0\qquad \text{for all } F\in\Mat,
\]
where $C_0$ and $C_1$ are positive constants independent
of $F\in\Mat$ where $p>3$ to ensure that deformations
with finite energy are uniformly continuous.

We can then denote the set of deformations of finite energy
by
\[
W^{\phi}=\{\, y\in C(\bar \Omega;\Vector) :
\int_{\Omega}\phi(\nabla y(x),\theta)\,dx<\infty\, \},
\]
and we can define the set $\A$ of admissible
deformations to be
\begin{equation}
\label{eq:admissible_deformations}
  \A =
 \{\, y \in W^\phi:
y(x) = y_0( x) \text{ for all } x \in \partial \Omega \,\}.
\end{equation}
Since we assume that the set  of admissible
deformations $\A$ is constrained on the entire boundary $\partial\Omega,$
we can neglect the loading energy in \eqref{energy}.  We will
also set the interfacial energy to be zero in this section so as to
consider the idealized model for which the length scale of the
microstructure is infinitesimally small.
For the theorems below, we assume boundary conditions compatible with a simple
laminate mixing
$QU_{i}$ for $Q\in\Rot$ with volume fraction $\lambda$ and $U_j$ with volume fraction
$1-\lambda,$
\[
 y_0(x)=\left[\lambda QU_{i} + (1-\lambda)U_{j}\right]x \qquad\text{for all } x\in\Omega,
\]
where for
$a\in\Real ^3$ and $n\in\Real^3,$ with $a,\,n\not =
0,$ we have the interface equation [1,\,19,\,25]
\begin{equation*}
\label{twin}
QU_i=U_j+a\otimes n.
\end{equation*}

We consider the finite element approximation of the variational problem
\[
\inf_{y\in\A}{\mathcal E}(y)
\]
given by
\[
\inf_{y_h\in\A_h}{\mathcal E}(y_h)
\]
where $\A_h$ is a finite-dimensional subspace of $\A$
defined for $h\in (0, h_0]$ for some $h_0>0.$
The following approximation theorem for the energy has been proven
 for the $P_k$ or $Q_k$ type conforming
finite elements on quasi-regular meshes, in particular
for the $P_1$ linear elements defined on
tetrahedra and the $Q_1$ trilinear elements defined on
rectangular parallelepipeds
[4,\,10,\,22,\,25--27].

{\bf Theorem 2.1.} \it
For each $h\in (0, h_0]$, there exists $y_h \in \A_h$ such that
\begin{equation}
\label{approx}
{\mathcal E}(y_h)=\min_{z_h\in \A_h} {\mathcal E} (z_h)\le Ch^{1/2}.
\end{equation}
\rm

We next define the volume fraction that an admissible deformation $y \in \A$
is in the $k$-th variant
$\Rot U_k$ for $k \in \{ 1, \, \ldots, N \}$ by
\[
  \tau_{k}(y) = \frac{\meas  \Omega_{k}(y)}{\meas  \Omega}
\]
where
\[
  \Omega_{k}(y) = \{\, x \in \Omega: \pi(\gy(x)) \in \Rot U_{k} \,\}.
\]

The following stability theory was first proven for the orthorhombic
to mono\-clinic transformation ($N=2$) [26]
and then for the cubic to tetragonal transformation ($N=3$) [22].
The analysis of stability is more difficult for larger $N$
since the additional
wells give the crystal more freedom to deform without the cost of additional energy.
In fact, for the tetragonal to monoclinic transformation ($N=4$) [6],
the orthorhombic to triclinic transformation ($N=4$) [16], and
the cubic to orthorhombic transformation ($N=6$) [4] we have shown that
there are special lattice constants for which the laminated microstructure is not stable.
Error estimates are obtained by substituting the approximation result
\eqref{approx} in the following stability results.

In each case for which we have proven the approximation of the microstructure
to be stable, we have derived the following basic stability estimate for the
approximation of a simple
laminate mixing
$QU_{i}$ and $U_j$ which bounds the volume fraction that $y\in\A$
is in the variants $k\not = i,j$
\begin{equation}
\label{stab}
\tau_k(y)\leq C\left(\E(y)^{\frac12}+\E(y)\right) \qquad \text{for all } k\not =i,j
\text{ and }
y\in
\A.
\end{equation}
For the theorems that follow, we shall assume that the lattice
parameters are such that the estimate \eqref{stab} holds.

The following theorem gives estimates for the strong convergence of
the projection of the deformation gradient parallel to the laminates
(the projection of the deformation gradient transverse to the laminates
does not converge strongly [25]), the strong convergence of the deformation,
and the weak convergence of the deformation gradient.

{\bf Theorem 2.2.} \it
(1) For any $w\in \Real^3$ such that $w\cdot n=0$ and $|w|=1$, we
    have the estimate for the strong convergence of the projection of
the deformation gradient
\[
\int_\Omega |\left(\nabla y(x)-\nabla y_0(x)\right)w|^2\, dx \leq
C\left(\E(y)+\E(y)^{\frac12}\right)\qquad\text{for all } y\in \A.
\]

(2) We have the estimate for the strong convergence of the
deformation
\[
\int_\Omega |y(x)-y_0(x)|^2\,dx \leq C\left(\E(y)+\E(y)^{\frac12}\right)\qquad\text{for
all } y\in \A.
\]

(3) For any Lipshitz domain $\omega\subset\Omega$, there exists a
    constant $C=C(\omega)>0$ such that we have the estimate for
the weak convergence of the deformation gradient
\[
\left\| \int_\omega (\nabla y(x)- \nabla y_0(x)) \,dx \right\| \leq
C\left(\E(y)^{\frac18}+\E(y)^{\frac12}\right)\qquad\text{for all } y\in \A.
\]
\rm

For fixed $i,j$ with $i\not = j$ we define a projection operator
$\pi_{ij}:\Mat\rightarrow \Rot U_i\cup \Rot U_j$ by
\[
\|F-\pi_{ij}(F)\|=\left\{ \|F-G\|:G\in \Rot U_i\cup \Rot U_j\right\}
\qquad\text{for all
} F\in
\Mat,
\]
and the operators $\Theta:\Mat\rightarrow
SO(3)$
and $\Pi:\Mat\rightarrow \{QU_i, U_j\}$ by the
unique decomposition
\[
\pi_{ij}(F)=\Theta(F)\Pi(F)\qquad\text{for all }  F\in
\Mat.
\]

The next theorem
shows that the deformation gradients of energy-minimizing
sequences  must oscillate between $QU_i$ and $U_j$.

{\bf Theorem 2.3.} \it
We have for all $y\in \A$ that
\[
\int_\Omega \|\nabla y(x)-\Pi(\nabla y(x))\|^2\, dx\leq
C\left(\E(y)+\E(y)^{\frac12}\right).
\]
\rm

We now present an estimate for the local volume fraction
that a deformation $y\in\A$ is near  $QU_i$ or
$U_j$. To describe this, we define the sets
\[
\omega_\rho^i(y)=\{\,x\in\omega:\Pi(\nabla y(x))=QU_i\text{ and }  \|\nabla
y(x)-QU_i\|\leq \rho\,\},
\]
\[
\omega_\rho^j(y)=\{\,x\in\omega:\Pi(\nabla y(x))=U_j\text{ and } \|\nabla
y(x)-U_j\|\leq \rho\,\},
\]
for any subset $\omega\in\Omega$, $\rho>0$, and $y\in \A$.
The next theorem demonstrates that the deformation gradients of
energy-minimizing sequences  must oscillate
with local volume fraction $\lambda$ near $QU_i$ and local volume
fraction $1-\lambda$ near $U_j$.

{\bf Theorem 2.4.} \it
For any Lipshitz domain $\omega\subset\Omega$ and for any $\rho>0$,
there exists a constant $C=C(\omega,\rho)>0$ such that for all $y\in
\A$
\[
\left|{\frac{\meas \omega_\rho^i(y)}{ \meas \omega}} -\lambda\right|+
\left|{\frac{\meas \omega_\rho^j(y)}{ \meas \omega}} -(1-\lambda)\right|\leq
C\left(\E(y)^{\frac18}+\E(y)^{\frac12}\right).
\]
\rm

We next give an estimate for the
 weak stability of nonlinear
functions of deformation gradients.

{\bf Theorem 2.5.} \it
We have for all $f:\Omega\times \Mat\rightarrow \Real$ and $y\in \A$ that
\[
\int_\Omega \{f(x,\nabla y(x))-[\lambda
 f(x,QU_i)+(1-\lambda)f(x,U_j)]\}\,dx
\leq
 C\|f\|_{\V}\left[\E(y)^{\frac14}+\E(y)^{\frac12}\right]
\]
where
\[
\|f\|_{\V}^2=\int_\Omega \left\{ (\esssup \|\nabla_F f(x,F)\|)^2+|\nabla
z_f(x)n|^2+z_f(x)^2\right\}\, dx< \infty
\]
with $z_f:\Omega\rightarrow \Real$ defined by
\[
z_f(x)=f(x,QU_i)-f(x,U_j)\qquad\text{for all } x\in\Omega.
\]
\rm

\section{A computational model for martensitic phase transformation}
 \label{section 3}\setzero
\vskip-5mm \hspace{5mm }

We have developed a computational model for the
quasi-static evolution of the
martensitic phase transformation of a single crystal thin film
[8].
Our thin film model [7] includes surface energy, as well as
 sharp phase
boundaries with finite energy.   The model also includes the nucleation of regions
of the high temperature phase (austenite) as the film is heated through
the transformation temperature and nucleation of regions of
the low temperature phase (martensite) as the film is cooled.
The nucleation step in our
algorithm is needed since the film would otherwise not transform.

For our
total-variation surface energy model, the bulk energy for a film of
thickness $h>0$ with reference configuration
$\Omega_h\equiv\Omega\times(-h/2,\,h/2)$, where $\Omega\subset\Real^2$
is a domain with a Lipschitz continuous
boundary $\partial\Omega,$ is given by the sum of the surface energy
and the elastic energy
\begin{equation}
\label{bulk}
  \kappa \int_{\Omega_h}|D(\nabla u)|
  +
  \int_{\Omega_h}\phi(\nabla u,\theta)\,dx,
\end{equation}
where $\int_{\Omega_h}|D(\nabla u)|$ is
the total variation of the deformation
gradient~[7] and
$\kappa$ is a small positive constant.

We have shown
in~[7]
that energy-minimizing deformations $u$ of the
bulk energy \eqref{bulk} are asymptotically of the form
\begin{equation*}
  u(x_1,x_2,x_3)
  =
  y(x_1,x_2)+b(x_1,x_2)x_3
  +
  \text{o}(x_3^2)
  \text{ for }(x_1,x_2)\in\Omega,\ x_3\in(-h/2,\,h/2),
\end{equation*}
(which is similar to that found for a diffuse interface model [3]) where
$(y,b)$ minimizes the thin film energy
\begin{equation}
\label{eq:energy_kappa_BV}
  \E(y,b,\theta)
  =
  \kappa\left(\int_\Omega|D(\nabla y|b|b)|
    +
    \sqrt{2}\int_{\partial\Omega}|b-b_0|\right)
  +
  \int_\Omega\phi(\nabla y|b,\theta)\,dx
\end{equation}
over all deformations
of finite energy such
that $y=y_0$ on $\partial\Omega.$
 The map $b$
describes the deformation of the cross-section
relative to the film~[3]. We denote by
$(\nabla y|b)\in\Mat$ the matrix whose first two columns are given by
the columns of $\nabla y$ and the last column by $b.$
In the above equation, $\int_\Omega|D(\nabla y|b|b)|$ is the total
variation of the vector valued function $(\nabla
y|b|b):\Omega\to\Real^{3\times4}.$

We describe our finite element approximation
of~(\ref{eq:energy_kappa_BV}) by letting the elements of a triangulation
$\tau$ of $\Omega$ be denoted by $K$ and the inter-element edges by
$e$. We denote the internal edges by $e\subset\Omega$ and
the boundary edges by $e\subset\partial\Omega.$
 We define the jump of a function $\psi$
across an internal edge
$e\subset\Omega$ shared by two elements
$K_1,K_2\in\tau$ to be
\[
  \jump{\psi}_e=\psi_{e,K_1}-\psi_{e,K_2},
\]
where $\psi_{e,K_i}$ denotes the trace on $e$ of $\psi|_{K_i},$
 and we define $\psi|_e$
to be the trace on $e$ for a boundary edge $e\subset\partial\Omega.$
 Next, we denote
by $\mathcal{P}_1(\tau)$ the space of continuous, piecewise linear
functions on $\Omega$ which are linear on each $K\in\tau$ and by
$\mathcal{P}_0(\tau)$ the space of piecewise constant functions on
$\Omega$ which are constant on each $K\in\tau$. Finally, for deformations
$(y,b)\in\mathcal{P}_1(\tau)\times\mathcal{P}_0(\tau)$ and
temperature fields $\tilde\theta\in\mathcal{P}_0(\tau)$, the
energy~(\ref{eq:energy_kappa_BV}) is well-defined and we have that
\begin{align*}
  &\kappa\left[\int_\Omega|D(\nabla y|b|b)|+
    \sqrt{2}\int_{\partial\Omega}|b-b_0|\right]
  +\int_\Omega\phi(\nabla y|b,\tilde\theta)\,dx
  \\
  &\qquad=\kappa\left(\sum_{e\subset\Omega}\Bigl|\jump{(\nabla y|b|b)}_e\Bigr|\,|e|
  +
  \sqrt{2}\sum_{e\subset\partial\Omega}\Bigl|b|_e-b_0|_e\Bigr|\,|e|\right)
  +
  \sum_{K\in\tau}\phi\big((\nabla y|b,\tilde\theta)|_K\big)\,|K|,
\end{align*}
where $|\cdot|$ denotes the
euclidean vector norm, $|e|$ denotes the length of the edge $e,$
  $|K|$ is the area of
the element $K,$ and
\[
  \Bigl|\jump{(\nabla y|b|b)}_e\Bigr|
  =
  \left(\bigl|\jump{\nabla y}_e\bigr|^2
  +
  2\,\bigl|\jump{b}_e\bigl|^2\right)^{1/2}.
\]
The above term is not differentiable everywhere, so we have
regularized it in our numerical simulations.

Since martensitic alloys are known to transform
on a fast time scale, we model the transformation of the film from martensite to
austenite  during heating by assuming that the film reaches an elastic
equilibrium on a faster time scale than the evolution of the
temperature, so the temperature
$\tilde\theta(x,t)$ can be obtained from
a time-dependent model for thermal evolution [8].
To compute the evolution of the deformation, we partition the time interval
$[0,T]$ for
$T>0$ by
$0=t_0<t_1<\dots<t_{L-1}<t_L=T$ and then obtain the solution
$(y(t_\ell),b(t_\ell))\in\mathcal A_\tau$ for $\ell=0,\dots,L$ by
computing a local minimum for the energy $\E(v,c,\theta(t_\ell))$
with respect to the space of approximate admissible deformations
\begin{equation}
\label{eq:A_0a}
  \mathcal A_\tau
  =
  \{(v,c)\in\mathcal{P}_1(\tau)\times\mathcal{P}_0(\tau):
    \ v=y_0\text{ on }\partial\Omega\}.
\end{equation}
Since the martensitic transformation strains $\U\subset\Mat$ are local
minimizers of the energy density $\phi(F,\theta)$ for all
$\theta$ near $\theta_T,$ a deformation that is in the martensitic phase will continue to
be a local minimum for the bulk energy $\E(v,c,\theta(t))$ for
$\theta>\theta_T.$  Hence, our computational
model will not simulate a transforming film if we compute
$(y(t_\ell),b(t_\ell))\in\mathcal A_\tau$ by using an energy-decreasing algorithm with
the initial state for the iteration at $t_\ell$ given by the deformation
at $t_{\ell-1}$, that is, if
$(y^{[0]}(t_{\ell}),b^{[0]}(t_{\ell}))=(y(t_{\ell-1}),b(t_{\ell-1}))$.
We have thus developed and utilized an algorithm to
nucleate regions of austenite into
$(y(t_{\ell-1}),b(t_{\ell-1}))\in\mathcal A_\tau$ to obtain an initial
iterate $(y^{[0]}(t_{\ell}),b^{[0]}(t_{\ell}))\in\mathcal A_\tau$ for
the computation of $(y(t_\ell),b(t_\ell))\in\mathcal A_\tau$.

We used an ``equilibrium distribution'' function,
$P(\theta)$, to determine the probability for which the crystal will be in
the austenitic phase at temperature $\theta$
and we assume that an equilibrium distribution has been reached
during the time between $t_{\ell-1}$ and $t_\ell$. The distribution
function $P(\theta)$ has the property that $0<P(\theta)<1$ and
\begin{gather*}
  P(\theta)\to0
  \text{ as }\theta\to-\infty
  \quad\text{ and }\quad
  P(\theta)\to1
  \text{ as }\theta\to\infty.
\end{gather*}

At each time $t_{\ell},$ we first compute a pseudo-random number $\sigma(K,\ell)\in (0,1)$ on
every triangle $K\in\tau,$ and we then compute
$(y^{[0]}(t_{\ell}),b^{[0]}(t_{\ell}))\in\mathcal A_\tau$ by
 ($x_K$ denotes the barycenter of $K$):

\begin{enumerate}
\item If $\sigma(K,\ell)\le P\left(\theta(x_K,t_{\ell})\right)$ and
  $\left(\nabla
    y(x_K,t_{\ell-1})|b(x_K,t_{\ell-1}),\theta(x_K,t_\ell)\right)$ is
  in \newline austenite, then set
\[
 (y^{[0]}(t_{\ell}),b^{[0]}(t_{\ell}))
=(y(t_{\ell-1}),b(t_{\ell-1}))
 \text{ on
}K.
\]
\item If $\sigma(K,\ell)\le P\left(\theta(x_K,t_{\ell})\right)$ and
  $\left(\nabla
    y(x_K,t_{\ell-1})|b(x_K,t_{\ell-1}),\theta(x_K,t_\ell)\right)$ is
  in \newline martensite, then transform to austenite on $K.$
\item If $\sigma(K,\ell)>P\left(\theta(x_K,t_{\ell})\right)$ and
  $\left(\nabla
    y(x_K,t_{\ell-1})|b(x_K,t_{\ell-1}),\theta(x_K,t_\ell)\right)$ is
  in \newline austenite, then transform to martensite on $K.$
\item If $\sigma(K,\ell)>P\left(\theta(x_K,t_{\ell})\right)$ and
  $\left(\nabla
    y(x_K,t_{\ell-1})|b(x_K,t_{\ell-1}),\theta(x_K,t_\ell)\right)$ is
  in \newline martensite, then set
\[
(y^{[0]}(t_{\ell}),b^{[0]}(t_{\ell}))
=(y(t_{\ell-1}),b(t_{\ell-1}))
 \text{ on
}K.
\]
\end{enumerate}

We have shown in [8] for a thin film of a
CuAlNi alloy in the ``tent'' configuration that we can
compute the nucleation above by
setting $y^{[0]}(t_{\ell})
  =
  y(t_{\ell-1})\in\mathcal P_1(\tau)$ and by updating the piecewise
constant
 $b^{[0]}(t_\ell)\in\mathcal
P_0(\tau)$
by
\begin{equation*}
  b^{[0]}(x_K,t_{\ell})
  =
  \dfrac{y_{,1}(x_K,t_{\ell-1})\times y_{,2}(x_K,t_{\ell-1})}
  {|y_{,1}(x_K,t_{\ell-1})\times y_{,2}(x_K,t_{\ell-1})|}
  \qquad\text{ on }K
\end{equation*}
to nucleate austenite and
\begin{equation*}
  b^{[0]}(x_K,t_{\ell})
  =
  \gamma\,\dfrac{y_{,1}(x_K,t_{\ell-1})\times y_{,2}(x_K,t_{\ell-1})}
  {|y_{,1}(x_K,t_{\ell-1})\times y_{,2}(x_K,t_{\ell-1})|}
  \qquad\text{ on }K
\end{equation*}
to nucleate martensite.

We then compute $(y(t_\ell),b(t_\ell))\in\mathcal A_\tau$ by the
Polak-Ribi\`ere conjugate gradient method with initial iterate
$(y^{[0]}(t_{\ell}),b^{[0]}(t_{\ell}))\in\mathcal A_\tau$.
We have also experimented with several other versions of the above
algorithm for the computation of $b^{[0]}(t_{\ell}).$ For example,
the above algorithm can be modified to utilize different probability
functions $P(\theta)$ in elements with increasing and decreasing
temperature. We can also prohibit the transformation from austenite
to martensite in an element in which the temperature is increasing or
prohibit the transformation from martensite to austenite in an
element for which the temperature is decreasing.

\label{lastpage}

\end{document}